\documentclass[12pt]{amsart}
\usepackage{amsmath}
\usepackage{amssymb}
\def\r{\mathrm}

\def\e{\varepsilon}

\def\Z{\mathbb Z}

\def\eqdef{{\buildrel \r{def} \over =}}

\def\oo{\infty}

\def\half{{\tfrac12}}

\def\upref#1{\upn{\ref{#1}}}

\theoremstyle{definition}

\theoremstyle{plain}
\newtheorem{thm}{Theorem}
\newtheorem{cor}{Corollary}
\newtheorem{lem}{Lemma}
\newtheorem{prp}{Proposition}
\theoremstyle{remark}
\newtheorem{rmrk}{Remark}
\newtheorem{ex}{Example}

\title[A general Hsu-Robbins-Erd\"os type estimate]{A general
Hsu-Robbins-Erd\"os type estimate of tail probabilities
of sums of independent identically distributed random variables}
\author{Alexander R. Pruss}
\date{April 1, 2002}
\address{Department of Philosophy\\
Georgetown University\\
Washington, DC 20057\\
U.S.A.}
\email{ap85@georgetown.edu}
\begin{document}
\begin{abstract}
	Let $X_1,X_2,\dots$ be a sequence of independent and identically
	distributed random variables, and put $S_n=X_1+\dots+X_n$.
	Under some conditions on the positive sequence $\tau_n$ and the
	positive increasing sequence $a_n$, we
	give necessary and sufficient conditions for the convergence of
$\sum_{n=1}^\oo \tau_n P(|S_n|\ge \e a_n)$
	for all $\e>0$, generalizing Baum and Katz's~(1965)
	generalization of the Hsu-Robbins-Erd\H os (1947, 1949) law of
	large numbers, also allowing us to characterize the
	convergence of the above series in the case where
	$\tau_n=n^{-1}$ and $a_n=(n\log n)^{1/2}$ for $n\ge 2$, thereby
	answering a question of Sp\u ataru.  Moreover, some results for
	non-identically distributed independent random variables are
	obtained by a recent comparison inequality. Our basic method is to use a
	central limit theorem estimate of Nagaev~(1965) combined with
	the Hoffman-J\o rgensen inequality~(1974).

\end{abstract}

\maketitle
\bibliographystyle{amsplain}
\baselineskip=1.92\baselineskip
\section{Introduction and main result}
Hsu and Robbins~\cite{HsuRobbins} showed that if
$X_1,X_2,\dotsc$ are independent and identically distributed mean zero random variables with
finite variance, then
\begin{equation}\label{eq:HR-concl}
	\sum_{n=1}^\infty P(|S_n| \ge \e n) < \infty, \qquad
	\forall \e >0,
 \end{equation}
where $S_n=X_1+\dots+X_n$.
Erd\H{o}s~\cite{Erdos,Erdos:Two} then showed
that a converse implication also holds, namely that if
\eqref{eq:HR-concl} holds  and the $\{ X_i \}$ are
independent and identically distributed, then $E[X_1^2]<\infty$ and $E[X_1]=0$.
Since then, a number of extensions in several directions have been
proved; see \cite{LRJW} and \cite{Pruss:Spataru} for partial
bibliographies and brief discussions.

The purpose of the present paper
is to prove a new and very general extension of the Hsu-Robbins-Erd\H os
law of large numbers.  Among other results, this extension will allow one
to give necessary and sufficient conditions (see Corollary~\ref{cor:Sp}
in Section~\ref{sec:cor})
for the truth of
\[
	\sum_{n=1}^\infty n^{-1} P(|S_n| \ge \e (n\log n)^{1/2}) < \infty, \qquad
	\forall \e >0,
 \]
where $S_n=X_1+\dots+X_n$ for independent and identically distributed
$X_1,X_2,\dots$, thereby answering a question of Sp\u ataru.

Throughout, terms like ``positive'' and ``increasing''
indicate the non-strict varieties (``non-negative'' and
``non-decreasing'', respectively).

To state our general result, we need a definition:
A positive sequence $\{ \tau_n \}_{n=1}^\oo$ is said to {\em satisfy
Condition~A} provided that for all positive decreasing
sequences $\{ c_n \}_{n=1}^\oo$
such that $\sum_{n=1}^\oo \tau_n
\min(nc_n,1)$ converges, we likewise have $\sum_{n=1}^\oo \tau_n nc_n$
converging.  The following Proposition gives us a simple sufficient criterion
for Condition~A.  A proof will be given in
Section~\ref{sec:proofs}, below.

\begin{prp}\label{prp:condA} Suppose that either
$\liminf_{n\to\oo}\tau_n>0$ or that both of the following two conditions
are satisfied:
\begin{enumerate}
\item[\upn{(a)}] $\liminf_{n\to\oo}n\tau_n > 0$, and
\item[\upn{(b)}] there is a constant $C\in (0,\oo)$ such that for each
$n\in\Z^+_0$, if\/ $2^{n-1}\le k<2^n$, then
$$
	C\tau_{2^{n-1}} \ge \tau_k \ge C^{-1} \tau_{2^n}.
$$
\end{enumerate}
Then, $\{ \tau_n \}$ satisfies Condition~\upn{A}.
\end{prp}

We now state our main result.  Conditions \eqref{eq:aux-cond} and
\eqref{eq:aux-aux-cond} in the result below are technical conditions
that will be satisfied for a wide class of sequences $\{ \tau_n \}$ and
$\{ a_n \}$.  In particular, they will automatically hold in the case
where $a_n=L(n) n^\alpha$ and $\tau_n=K(n) n^\beta$, assuming $L$ and $K$
are slowly-varying functions, $\alpha>1/3$ and $\liminf n\tau_n>0$ (see
Proposition~\ref{prp:sv}, below).

\begin{thm}\label{th:main} Let $X_1,X_2,\dotsc$ be independent
identically distributed random variables.
Put $S_n=X_1+\dots+X_n$.
Let $\{ \tau_n \}_{n=1}^\oo$ be any positive sequence of numbers,
and let $\{ a_n \}_{n=1}^\oo$ be an increasing strictly positive
sequence tending to $\oo$.
Furthermore, assume there exist finite constants
$\theta\ge 1$, $C>0$ and $N\ge 2$ such that for all $n\ge N$ we have:
\begin{equation}\label{eq:aux-cond}
	\frac{a_n^{3\theta}}{n^{\theta-1}} \sum_{k=n}^\oo
	\frac{k^{\theta}\tau_k}{a_k^{3\theta}}
	  \le C \sum_{k=1}^{n-1} k\tau_k.
\end{equation}
If $\theta>1$, then additionally assume that
\begin{equation}\label{eq:aux-aux-cond}
	\liminf_{n\to\oo} \inf_{k\ge n} \frac{a_k^3}{k a_n^3}
	   \sum_{j=1}^{n-1} j\tau_j > 0.
\end{equation}

Then, if all of the following three conditions hold:
\begin{enumerate}
\item[\upn{(i)}]
	there is a sequence $\{ \mu_n \}_{n=1}^\oo$ with $\mu_n$ a
	median of $S_n$ for all $n$, such that for all $\e>0$,
        we have $\sum_{n\in M(\e)}
        \tau_n<\oo$, where $M(\e)=\{ n\in \Z^+ : |\mu_n| > \e a_n \}$,
\item[\upn{(ii)}] $\sum_{n=1}^\oo n\tau_n P(|X_1|\ge \e a_n)<\oo$ for
all $\e>0$, and
\item[\upn{(iii)}] $\sum_{n=1}^\oo \tau_n
e^{-\e^2 a_n^2 / (nT_{\e,n})} <\oo$ for all $\e>0$,
where $T_{\e,n}=E[X_1^2\cdot 1_{\{|X_1|<\e a_n\}}]$,
and where $e^{-t/0}=0$ for all $t>0$,
\end{enumerate}
we will have:
\begin{equation}\label{eq:conv}
	\sum_{n=1}^\oo \tau_n P(|S_n|\ge \e a_n) < \oo, \quad \forall
	\e>0,
\end{equation}
and, conversely, if the sequence $\{ \tau_n \}$ satisfies Condition~A and
\eqref{eq:conv} holds, conditions~\upn{(i)}, \upn{(ii)} and~\upn{(iii)}
will also hold.
\end{thm}

The proof will be given in Section~\ref{sec:proofs}, below.  The result
is closely related to work of Klesov~\cite[Theorem~4]{Klesov}, though we
are working with a more general class of sequences $\{ a_n \}$.  Our
proof will be based upon Klesov's Hoffman-J\o rgensen inequality based
approach, combined with a central limit theorem estimate of
Nagaev~\cite{Nagaev}, the latter being used rather like in
\cite{Pruss:Riemann}.

The particular newness of the result is that it works for $a_n$ near
the critical growth $n^{1/2}$ involved in the central limit theorem. For
instance, as already advertised, in Corollary~\ref{cor:Sp} we will use
the Theorem to characterize the cases where \eqref{eq:conv} holds with
$\tau_n=n^{-1}$ and $a_n = (n\log n)^{1/2}$ (for $n\ge 2$).

Now, recall that a measurable function $\phi$ on $[0,\oo)$ is slowly varying (in
the sense of Karamata) providing that for all $\lambda>0$ we have
$\lim_{x\to\oo} \phi(\lambda x)/\phi(x)=1$.

\begin{prp}\label{prp:sv} If $K$ and $L$ are strictly positive
slowly varying functions on $[0,\oo)$, and if $a_n=K(n) n^\alpha$ for some
$\alpha>\tfrac13$ and $\tau_n=L(n) n^\beta$ for some real $\beta$ such
that $\liminf_n n\tau_n > 0$, then
conditions \eqref{eq:aux-cond} and \eqref{eq:aux-aux-cond} will be satisfied
for a sufficiently large $\theta$.
\end{prp}

\begin{proof}[Sketch of proof of Proposition~\upref{prp:sv}]
First, for sufficiently large $\theta$, approximate $\sum_{k=n}^\oo
a_k^{-3\theta} k^\theta \tau_k$ by an integral and use the fact that if
$\phi$ is slowly varying then $\int_{x}^\oo t^{-p} \phi(t)\, dt$ is
asymptotic to $(p-1)^{-1} x^{1-p} \phi(x)$ for large $x$ if $p>1$ by
Karamata's integral theorem
\cite[Proposition~1.5.10]{BGT}. It is here that the condition that
$\alpha>\tfrac13$ will be used;  one will also need to use the easy fact
that powers and products of slowly varying functions are slowly varying.
This will show via a straightforward computation that the left
hand side of \eqref{eq:aux-cond} is asymptotic to $C n^{2+\beta} K(n)$.
But also using an integral to approximate the right hand side of
\eqref{eq:aux-cond} and using the fact that if $\phi$ is slowly varying
then $\int_1^x t^p \phi(t) \, dt$ is asymptotic to $(p+1)^{-1} x^{p+1}
\phi(x)$ for large $x$ by Karamata's other integral theorem~\cite[Proposition~1.5.11]{BGT}
if $p > -1$, we can see that the right hand side of \eqref{eq:aux-cond} is
also asymptotic to $C n^{2+\beta} K(n)$, and so \eqref{eq:aux-cond} holds.

Condition~\eqref{eq:aux-aux-cond} can be obtained by noting that $K$ can
without loss of generality be replaced by a normalized slowly varying
function (see \cite[p.~15]{BGT}) asymptotic to $K$, and then $a_k^3/k$
will be an increasing function for sufficiently large $k$ if $\alpha>
1/3$, by the Bojanic-Karamata theorem~\cite[Theorem~1.5.5]{BGT}, so that
the infimum on the left hand side of \eqref{eq:aux-aux-cond} will be
attained at $k=n$, and the truth of \eqref{eq:aux-aux-cond} will follow
from the condition that $\liminf n\tau_n > 0$.	Condition~A will also be
satisfied by $\{ \tau_n \}$ under the above circumstances, as is very
easy to see by Proposition~\ref{prp:condA}. \end{proof}

\begin{rmrk}\label{rk:MS}  By a maximal inequality of
Montgomery-Smith~\cite{MS}, condition~\eqref{eq:conv} is equivalent to: $$
\sum_{n=1}^\oo \tau_n P\bigl(\bigl|\sup_{1\le k\le n} S_k\bigr|\ge
\e a_n\bigr)<\oo,
\qquad\forall{\e>0}. $$
\end{rmrk}

\begin{rmrk}
Readers familiar with Hsu-Robbins-Erd\H os laws of large numbers may be
surprised at condition~(iii) in Theorem~\ref{th:main}, since normally
these laws of large numbers simply have \eqref{eq:conv} equivalent to
(i) and (ii) (under appropriate conditions on $\{ \tau_n \}$ and $\{ a_n \}$).
In general, condition~(iii) cannot be eliminated from
Theorem~\ref{th:main}, and seems to become particularly significant for
$a_n$ near the critical grown $n^{1/2}$.
(Example~\ref{ex:MS} in Section~\ref{sec:cor} will
show that condition~(iii) cannot be eliminated if $\tau_n=n^{-1}$ and
$a_n=(n\log n)^{1/2}$.)  However, in a number of special cases,
condition~(iii) can indeed be removed, as is seen in the following
result.
\end{rmrk}

\begin{thm}\label{th:aux}  Suppose the $\{ \tau_n \}$ is positive and $\{ a_n \}$ is increasing
and strictly positive, and that there are finite real constants
$\theta\ge 1$, $C>0$ and $N\ge 2$ such that for all $n\ge N$ we have:
\begin{equation}\label{eq:aux-cond2}
	\frac{a_n^{2\theta}}{n^{\theta-1}} \sum_{k=n}^\oo
	\frac{k^{\theta}\tau_k}{a_k^{2\theta}}
	  \le C \sum_{k=1}^{n-1} k\tau_k
\end{equation}
and
\begin{equation}\label{eq:aux-aux-cond2}
	\liminf_{n\to\oo} \inf_{k\ge n} \frac{a_k^2}{k a_n^2}
	   \sum_{j=1}^{n-1} j\tau_j > 0.
\end{equation}
Then, any random variable $X_1$ satisfying condition~\upn{(ii)} of
Theorem~\upref{th:main} automatically satisfies condition~\upn{(iii)} of
that Theorem.
\end{thm}

The proof will be given in Section~\ref{sec:proofs}, below.  Note that
\eqref{eq:aux-cond2} automatically implies \eqref{eq:aux-cond}, and
\eqref{eq:aux-aux-cond2} likewise implies \eqref{eq:aux-aux-cond}.

\begin{rmrk}\label{rk:sv-2}
It is not difficult to see that if $K$ and $L$ are slowly varying
strictly positive functions on $[0,\oo)$, and if $a_n=K(n) n^\alpha$ for
some $\alpha> \tfrac12$ while $\tau_n=L(n) n^\beta$ for some
real $\beta$ such that $\liminf_n n\tau_n > 0$, then conditions \eqref{eq:aux-cond2} and
\eqref{eq:aux-aux-cond2} will be satisfied for a sufficiently large
$\theta$.  Since by Proposition~\ref{prp:sv}, conditions~\eqref{eq:aux-cond} and
\eqref{eq:aux-aux-cond} would allow $\alpha$ to be $\half$ (or in fact any
value greater than $\tfrac13$), this helps to further illustrate how
condition~(iii) becomes relevant close to the critical growth $n^{1/2}$
of $a_n$,
but for faster growths (i.e., around $n^\alpha$ for $\alpha>\half$), it
can be eliminated by Theorem~\ref{th:aux}. \end{rmrk}

\section{Corollaries and applications}\label{sec:cor}
We now obtain the following generalization due to
Baum and Katz~\cite{BaumKatz} of the Hsu-Robbins-Erd\H os
result, thereby showing that
Theorem~\ref{th:main} is indeed more general than
the Hsu-Robbins-Erd\H os law of large numbers. (See~\cite{LRJW} for a
discussion of the pedigree of the Baum and
Katz result.)
\begin{cor}\label{cor:bk}
Suppose that $X_1,X_2,\dotsc$ are independent identically
distributed random variables.  Let $S_n=X_1+\dots+X_n$.
Fix $r\ge 1$ and $0<p<2$.
Then the conjunction of the conditions
\begin{enumerate}
\item[\upn{(a)}] if $p\ge 1$ then $E[X_1]=0$, and
\item[\upn{(b)}] $E[|X_1|^{rp}]<\oo$
\end{enumerate}
holds if and only if
\begin{equation}\label{eq:bk-concl}
	\sum_{n=1}^\oo n^{r-2} P(|S_n|\ge \e n^{1/p})<\oo, \qquad
	\forall{\e>0}.
\end{equation}
\end{cor}

\begin{proof}[Proof]
Let $\tau_n=n^{r-2}$ and $a_n=n^{1/p}$.
Observe that \eqref{eq:aux-cond}, \eqref{eq:aux-aux-cond},
\eqref{eq:aux-cond2}, \eqref{eq:aux-aux-cond2} and Condition~A all
hold for appropriate choices of
$\theta$ (one can use Proposition~\ref{prp:sv} and Remark~\ref{rk:sv-2} here if
one so desires, but in fact a direct verification is easy).
Then by the Marcinkiewicz-Zygmund strong law of large numbers~\cite[\S 16.4.A.3]{Loeve},
if (a) and (b)
hold, then $S_n/a_n\to 0$ almost surely, hence also in probability, and
therefore condition~(i) of
Theorem~\ref{th:main} holds as well.  Moreover, (ii) is equivalent to (b),
and (ii) implies (iii) by Theorem~\ref{th:aux}.
Thus, by Theorem~\ref{th:main}, if (a)
and (b) hold, \eqref{eq:bk-concl} follows.  Conversely, if
\eqref{eq:bk-concl} holds, then since $\tau_n$ satisfies
Condition~A by Proposition~\ref{prp:condA}, it follows from
Theorem~\ref{th:main} that (i) and (ii)
(and (iii), but that is not needed) hold.  Hence (b) holds, since it is
equivalent to (ii).  It remains to show that (a) holds.
The easiest way to do this is to note that if $p\ge 1$ then by
the Marcinkiewicz-Zygmund strong law of large numbers, we have
$(S_n-nE[X_1])/a_n\to 0$ almost surely, hence in probability, and
therefore $(\mu_n-nE[X_1])/n^{1/p}\to 0$, where
$\mu_n$ is a median of $S_n$, and by condition~(i) it will then follow
that $E[X_1]=0$.
\end{proof}

Recall that random variables $X_1,\dots,X_n$ are said to be $K$-{\em
weakly
mean dominated} by a random
variable $X$ providing that for all $\lambda$ we have:
$$
	\frac1n \sum_{k=1}^n P(|X_k|\ge \lambda) \le K P(|X|\ge\lambda)
$$
(see \cite{Gut}).

\begin{cor}\label{cor:wmd}
Let $X$ be a random variable and $K$ any finite constant.
Let $\{ X_{k,n} \}_{1\le k\le n;\, n\ge 1}$ be a triangular array of
random variables such that
$X_{1,n},\dots,X_{n,n}$ are independent and $K$-weakly mean dominated
by $X$ for each fixed $n$.
Let $\{ \tau_n \}$ be a positive sequence of
numbers and let $\{ a_n \}$ be an increasing strictly positive sequence for
$n\ge 1$.
Put $S_n=X_{1,n}+\dots+X_{n,n}$.
Suppose that condition~\upn{(i)} of
Theorem~\upref{th:main} holds, \eqref{eq:aux-cond} is satisfied and if $\theta>1$ then
\eqref{eq:aux-aux-cond} holds as well.	Assume that
conditions~\upn{(ii)} and~\upn{(iii)} of
Theorem~\upref{th:main} are satisfied with $X$ in place of $X_1$.  Then,
$$
	\sum_{n=1}^\oo \tau_n P(|S_n|\ge \e a_n) < \oo
$$
for all $\e>0$.
\end{cor}

\begin{proof}
By Theorem~\ref{th:main} we have
\begin{equation}\label{eq:primed-conv}
	\sum_{n=1}^\oo \tau_n P(|S_n'|\ge \e a_n)<\oo,
	\qquad\forall\e>0,
\end{equation}
where $S_n'$ is the sum of $n$ independent copies of $X$ (condition~(i)
in this case will hold trivially by symmetry, while (ii) and (iii) were
assumed in the statement of Corollary~\upref{cor:wmd}.) One may slightly
modify the comparison result in \cite[Corollary~1]{Pruss:noniid} by
assuming our Theorem~\ref{th:main}'s condition~(i) in place of the
assumption in that paper that $S_n/a_n\to 0$ in
probability, which modification only very slightly affects the proof
(one will need to use \eqref{eq:for-ref}, below, after obtaining the
convergence of \cite[series~(1.2)]{Pruss:noniid} in the original proof
in \cite{Pruss:noniid}).
Thus modified, \cite[Corollary~1]{Pruss:noniid} together with \eqref{eq:primed-conv}
yields the conclusion of our Corollary~\ref{cor:wmd}.
\end{proof}

The following Corollary yielding a result similar to one of Hu, Moricz and
Taylor~\cite{HMT} (cf.\ \cite{Gut,HSV,HV}) can be derived from
Corollary~\ref{cor:wmd} exactly in the way that Corollary~\ref{cor:bk}
was derived from Theorem~\ref{th:main}.
\begin{cor}
Let $X$ be a random variable and let $K$ be any finite constant.
Let $\{ X_{k,n} \}_{1\le k\le n;\, n\ge 1}$ be a triangular array of
random variables such that $X_{1,n},\dots,X_{n,n}$ are independent $K$-weakly mean dominated
by $X$.  Let $S_n=X_{1,n}+\dots+X_{n,n}$.
Fix $r\ge 1$ and $0<p<2$.
Suppose $E[|X|^{rp}]<\oo$.  If $p\ge 1$ then assume also that
$E[X_{n,1}+\dots+X_{n,n}]=0$ for all $n$.
Then,
$$
	\sum_{n=1}^\oo n^{r-2} P(|S_n|\ge \e n^{1/p}), \qquad
	\forall{\e>0}.
$$
\end{cor}

Now define $\log^+ x=\log (2+x)$.  It is easy to see that
Theorem~\ref{th:main} implies the following result.
\begin{cor}\label{cor:Sp}
Let $X_1,X_2,\dotsc$ be independent and identically distributed random
variables.  Then,
\begin{equation}\label{eq:sp-conv}
	\sum_{n=2}^\oo n^{-1} P(|S_n| \ge \e (n\log n)^{1/2}) < \oo,
	\quad \forall \e>0,
\end{equation}
if and only if all of the following three conditions hold:
\begin{enumerate}
\item[\upn{(a)}] $E[X_1]=0$,
\item[\upn{(b)}] $E[X_1^2/\log^+ |X_1|]<\oo$, and
\item[\upn{(c)}] $\sum_{n=2}^\oo n^{-1-\e^2/T_{\e,n}}<\oo$ for
all $\e>0$, where $T_{\e,n}=E[X_1^2\cdot
1_{\{|X_1|<\e (n\log n)^{1/2}\}}]$.
\end{enumerate}
\end{cor}

The proof of the following Lemma will be given at the end of
Section~\ref{sec:proofs}.
\begin{lem}\label{lem:Sp}
	If condition \upn{(b)} of Corollary~\upn{\ref{cor:Sp}}
	holds, then $(S_n-E[S_n]) / (n\log n)^{1/2} \to 0$ in probability
	as $n\to\oo$.
\end{lem}

\begin{rmrk}
It is not known whether ``in probability'' can be replaced by ``almost
surely'' in Lemma~\ref{lem:Sp}.
\end{rmrk}

\begin{proof}[Proof of Corollary~\upref{cor:Sp}]
Let $\tau_n=1/n$ and $a_n=(n\log n)^{1/2}$ for $n\ge 2$.
Note that \eqref{eq:aux-cond} is easily seen to be
satisfied with $\theta=1$, and Condition~A holds by
Proposition~\ref{prp:condA}.  Observe that (b) is
equivalent to condition~(ii) of Theorem~\ref{th:main} in the present
setting, and that if (b) holds,
then (a) is equivalent to (i) by
Lemma~\ref{lem:Sp}, so that the conjunction of conditions (a) and (b)
of Corollary~\ref{cor:Sp}
is
equivalent to that of conditions (i) and (ii) of Theorem~\ref{th:main}.
Also, it is easy to check that (c) is
equivalent to (iii), since the summands in the sums in both conditions
are equal.
\end{proof}

Professor Aurel Sp\u ataru has asked the author whether
\eqref{eq:sp-conv} is equivalent to the conjunction of (a) and (b). This
would be expected by analogy with other Hsu-Robbins-Erd\H os laws of
large numbers (such as Corollary~\ref{cor:bk}). It
is this question that has inspired the present paper. In light of the
Corollary~\ref{cor:Sp}, Sp\u ataru's question is equivalent to asking
whether the conjunction of (a) and (b) implies (c).  The following
counterexample that Professor Stephen Montgomery-Smith has privately
communicated to the author  shows that the answer is negative, and hence so is
the answer to Sp\u ataru's question.

\begin{ex}\label{ex:MS}
Let $\{ K_m \}_{m=0}^\oo$ be a very rapidly increasing strictly positive
sequence, with $K_0=0$.
The degree of rapidity of increase will be chosen later so as to
make the argument go. Let $\psi(t)=(t \log t)^{1/2}$ for $t\ge 2$.
Extend $\psi$ linearly to the interval $[0,2]$ in such a way that $\psi(0)=0$.
Let $\phi$
be the inverse function of $\psi$.  Assume that $K_1\ge 2$. Let
$X_1$ be a random variable such that
$P(X_1=\psi(K_m))=P(X_1=-\psi(K_m))=2^{-m-1}/K_m$ for all $m$, and with
$P(|X_1| \notin \{ 0 \} \cup \{ \psi(K_m):m\in\Z^+ \})=0$.

Let $\e=1$.  We have
$$
	E[\phi(|X_1|)]
	 =\sum_{m=1}^\oo 2^{-m} K_m^{-1} \cdot K_m = 1.
$$
It is easy to see that in general $E[\phi(|X_1|)]<\oo$ if and only if (b) holds,
and hence indeed (b) is satisfied in the present case.	So is (a), since $X_1$
is symmetric (i.e., $X_1$ and $-X_1$ have the same distribution)
 and so we can put $\mu_n=0$ for all $n$.
Now, let $M(n)=\max \{ m\in \Z^+_0 : K_m< n \}$.  Then, with $T_{1,n}$
as in condition~(iii) of Theorem~\ref{th:main}, and as $K_m\ge 2$ for
all $m\ge 2$
\[
\begin{split}
   T_{1,n}&=\sum_{m=1}^{M(n)} (2^{-m}/K_m) (\psi(K_m))^2 \\
   &= \sum_{m=1}^{M(n)} (2^{-m}/K_m) (K_m \log K_m) \\
   &= \sum_{m=1}^{M(n)} 2^{-m}\log K_m \ge 2^{-M(n)}\log K_{M(n)}.
\end{split}
\]
Note that if $K_m <n \le K_{m+1}$, then $M(n)=m$ and so $T_{1,n}\ge
2^{-m}\log K_m$.  Thus, for $m\ge 1$ we have
\begin{equation}\label{eq:oneterm}
	\sum_{n=K_m+1}^{K_{m+1}}
	n^{-1-1/T_{1,n}}
	\ge\sum_{n=K_m+1}^{K_{m+1}}
	 n^{-1-2^m/\log K_m}.
\end{equation}
Now, for any $K\ge 2$ and $m\in\Z^+$, let $L(K,m)$ be an integer greater than $K$ and
sufficiently large that:
\begin{equation}
	\sum_{n=K+1}^{L(K,m)}
	  n^{-1-2^m/\log K}
	  \ge \frac12 \sum_{n=K+1}^\oo
	   n^{-1-2^m/\log K}.
\end{equation}
Such an $L(K,m)$ exists because the sum on the right hand side
converges.  Observe furthermore that
\begin{equation}\label{eq:threeterm}
\begin{split}
	\sum_{n=K+1}^\oo
	   n^{-1-2^m/\log K} \ge C \cdot 2^{-m} (\log K) (K+1)^{-2^m/\log K},
\end{split}
\end{equation}
for an absolute constant $C>0$ independent of $K\ge 2$ and $m\ge 0$. Combining
\eqref{eq:oneterm}--\eqref{eq:threeterm} we see that if $K_{m+1}\ge
L(K_m,m)$, then we have
\begin{equation}\label{eq:the-ineq00}
\begin{split}
	\sum_{n=K_m+1}^{K_{m+1}}
	n^{-1-1/T_{1,n}}
	 &\ge \frac C2 \cdot 2^{-m} (\log K_m) (K_m+1)^{-2^m/\log K_m}
	 \\
	 &= 2^{-m-1} C (\log K_m) \exp\left(-2^m
\frac{\log(K_m+1)}{\log K_m}\right).
\end{split}
\end{equation}
Inductively choosing the $K_m$ so that for all $m$
we have both
$$
2^{-m-1} (\log K_m) \cdot \exp\left(-2^m
\frac{\log(K_m+1)}{\log K_m}\right) \ge 1,
$$
and $K_{m+1}\ge
L(K_m,m)$, we then find by \eqref{eq:the-ineq00} that
that
\[
	\sum_{n=2}^\oo
	 n^{-1-1/T_{1,n}}
	\ge \sum_{m=1}^\oo \sum_{n=K_m+1}^{K_{m+1}}
	n^{-1-1/T_{1,n}}
	\ge C \sum_{m=1}^\oo 1 = \oo.
\]
Hence Corollary~\ref{cor:Sp}'s condition~(c) fails, and so we do have our
desired counterexample satisfying
(a), (b) but not (c), and hence by that Corollary, with inequality
\eqref{eq:sp-conv} also failing.
\end{ex}

Although the answer to Sp\u ataru's question is negative, we do have the
result under a slightly stronger moment condition than Corollary~\ref{cor:Sp}'s
condition~(b).

\begin{cor}\label{cor:Sp-weak}
Let $X$ be a random variable and let $K$ be any finite constant.
Let $\{ X_{k,n} \}_{1\le k\le n;\, n\ge 1}$ be a triangular array of
random variables such that $X_{1,n},\dots,X_{n,n}$ are $K$-weakly mean dominated
by $X$.
Assume that
\begin{enumerate}
\item[\upn{(a)}] $E[X_{n,1}+\dots+X_{n,n}]=0$ for all $n$, and
\item[\upn{(b)}] $E\left[\frac{X^2 (\log^+\log^+ |X|)^{1+\delta}}{\log^+ |X|}\right]<\oo$
for some $\delta>0$.
\end{enumerate}
Then,
$$
	\sum_{n=2}^\oo \frac1n P(|S_n| \ge \e (n\log n)^{1/2}) < \oo,
	\quad \forall \e>0.
$$
\end{cor}

\begin{proof} Using Corollary~\ref{cor:wmd} and the same methods as in
Corollary~\ref{cor:Sp} it suffices to show that if condition~(b)
of Corollary~\ref{cor:Sp-weak} holds, then condition~(c) of
Corollary~\ref{cor:Sp} is satisfied with $X$ in place of $X_1$.
To show this, without loss of generality (rescaling $X$ if necessary)
assume $\e=1$.
Note that if $T_{1,n}=E[X^2\cdot 1_{\{ |X|<a_n
\}}]$, where $a_n=(n\log n)^{1/2}$, then
\begin{equation}\label{eq:star2}
\begin{split}
	T_{1,n} &\le E\left[\frac{X^2 (\log^+\log^+ |X|)^{1+\delta}}{\log^+
	       |X|} \right]
		\cdot \frac{\log^+ a_n}{(\log^+\log^+ a_n)^{1+\delta}} \\
		&\le C \frac{\log^+ a_n}{(\log^+\log^+ a_n)^{1+\delta}} \\
		&\le C' \frac{\log^+ n}{(\log^+ \log^+ n)^{1+\delta}},
\end{split}
\end{equation}
where $\delta$ is as in (b), while $C$ and $C'$ are
strictly positive finite constants independent of $n\ge 2$ (but
dependent on the value of the expectation in (b)).  Hence,
if $N\ge 2$ is sufficiently large that $(\log^+\log^+ n)^{1+\delta}\ge 2C'
\log\log n$ for all $n\ge N$, then we have:
\[
\begin{split}
	\sum_{n=N}^\oo n^{-1-1/T_{1,n}}
	 &=\sum_{n=N}^\oo
	    n^{-1} \exp(-(\log n)/T_{1,n}) \\
	 &\le \sum_{n=N}^\oo
	    n^{-1} \exp(-(\log^+ \log^+
	    n)^{1+\delta}/C') \\
	 &\le \sum_{n=N}^\oo n^{-1}
	  \exp(-2\log\log n)
	 = \sum_{n=N}^\oo \frac1{n (\log n)^{2}} < \oo,
\end{split}
\]
by \eqref{eq:star2} and the choice of $N$.
Thus condition~(c) of Corollary~\ref{cor:Sp} is indeed satisfied
with $X$ in place of $X_1$.
\end{proof}

\section{Proofs and auxiliary results}\label{sec:proofs}
\begin{proof}[Proof of Proposition~\upref{prp:condA}]
If $\liminf \tau_n>0$ and $\sum_{n=1}^\oo \tau_n
\min(1,nc_n)<\oo$, then there are only finitely many $n\in\Z^+$
for which $nc_n\ge 1$, and hence $\sum_{n=1}^\oo \tau_n nc_n$ must also
converge since it differs from $\sum_{n=1}^\oo \tau_n \min(1,nc_n)$ only
in finitely many terms.

It remains to show that if (a) and (b) hold, then Condition~A.
To do this, suppose $c_n$ is a decreasing sequence such that
\begin{equation}\label{eq:c0}
\sum_{n=1}^\oo \tau_n\min(1,nc_n)<\oo.
\end{equation}
Then, using (b) we have:
\begin{equation}\label{eq:star00}
\begin{split}
	\oo& >\sum_{n=1}^\oo \tau_n\min(1,nc_n) \\
	   &\ge \sum_{k=0}^\oo C^{-1}2^k\tau_{2^{k+1}}\min(1,2^k
	   c_{2^{k+1}}) \\
	   &\ge \frac{1}{2C}\sum_{k=0}^\oo
	   2^{k}\tau_{2^{k+1}}\min(1,2^{k+1} c_{2^{k+1}}).
\end{split}
\end{equation}
Now, $\liminf 2^k \tau_{2^{k+1}}>0$ by (a), and hence it follows that only
finitely many of the $\min(1,2^{k+1} c_{2^{k+1}})$ can equal $1$
(since otherwise the right hand side of \eqref{eq:star00} would be
infinite), so
that except for at most finitely many values of $k$, we have $\min(1,2^{k+1}
c_{2^{k+1}})=2^{k+1}c_{2^{k+1}}$.  It thus
follows from \eqref{eq:star00} and (b)
that
\[
\begin{split}
	\oo & > \frac12 \sum_{k=0}^\oo 2^{k+1} \tau_{2^{k+1}} \cdot
	2^{k+1}c_{2^{k+1}} \\
	    &\ge \frac{1}{2C}\sum_{k=0}^\oo \sum_{j=0}^{2^{k+1}-1}
	      \tau_{j+2^{k+1}} {\frac{j+2^{k+1}}{2}}\cdot c_{j+2^{k+1}} \\
	    &= \frac{1}{4C}\sum_{n=2}^\oo \tau_n n c_n,
\end{split}
\]
and the proof is complete.
\end{proof}

\begin{lem}\label{lem:1}  In the setting of Theorem~\upref{th:main}, if the series
$\{ \tau_n \}$ satisfies Condition~A, then condition \eqref{eq:conv}
entails \upn{(i)} and \upn{(ii)}.
\end{lem}

\begin{proof}
Assume \eqref{eq:conv} holds.
Fix $\e>0$ and any
sequence of medians $\mu_n$ of the $S_n$.  Let $M$ be as in
condition~(i).	Then, it is easy to see that for $n\in M$ we have
$P(|S_n|\ge \e a_n)\ge 1/2$.  The convergence of $\sum_{n\in M} \tau_n$
follows immediately from this and \eqref{eq:conv}, and so condition~(i)
holds.

On the other hand, by \eqref{eq:conv} and the remark following the main theorem
in~\cite{Pruss:Klesov}, we have
$$
	\sum_{n=1}^\oo \tau_n \min (1,nP(|X_1|\ge \e a_n)) < \oo
$$
for all $\e>0$.  Condition~(ii) follows immediately from this together
with the fact that
$\{ \tau_n \}$ satisfies Condition~A while $\{ a_n \}$ is increasing.
\end{proof}

\begin{lem}\label{lem:elementary}
Suppose that $\{ \tau_n \}$ and $\{ \rho_n \}$ are positive  and that
$\{ b_n \}$ is increasing and strictly positive.  Let $X$ be any
random variable. If there is a constant $C\in (0,\oo)$
such that for all $n\ge 2$:
\begin{equation}\label{eq:aux-cond-gen}
	b_n^t \sum_{k=n}^\oo \rho_k  \le C \sum_{k=1}^{n-1} k\tau_k,
\end{equation}
then
$$
	\sum_{n=2}^\oo \rho_n E[|X|^t\cdot 1_{\{ |X|<b_n \}}]
\le C	     \sum_{n=1}^\oo n\tau_n P(|X|\ge b_n).
$$
\end{lem}

\begin{proof}[Proof of Lemma~\upref{lem:elementary}]
Assume without loss of generality that $b_1=0$.
Then, by Fubini's theorem and \eqref{eq:aux-cond-gen}:
\[
\begin{split}
  \sum_{n=2}^\oo \rho_n E[|X|^t\cdot 1_{\{ |X|<b_n \}}]
   &\le \sum_{n=2}^\oo \rho_n \sum_{k=1}^{n} b_k^t P(b_{k-1}\le |X|<b_k) \\
   &=\sum_{k=2}^\oo P(b_{k-1}\le |X|<b_k) b_k^t \sum_{n=k}^\oo \rho_n \\
   &\le C\sum_{k=2}^\oo P(b_{k-1}\le |X|<b_k) \sum_{n=1}^{k-1} n\tau_n \\
   &=C \sum_{n=1}^\oo n\tau_n \sum_{k=n+1}^\oo P(b_{k-1}\le |X|<b_k) \\
   &=C \sum_{n=1}^\oo n\tau_n P(|X|\ge b_n).
\end{split}
\]
\end{proof}

\begin{lem}\label{lem:comp} Suppose that $\{ \alpha_n \}$ and $\{ \beta_n \}$ are sequences
in $[0,1]$, that $\{ \tau_n \}$ is a positive sequence, and
that $r\in\Z^+$ is such that\/ $\sum_{n=1}^\oo \tau_n
|\alpha_n-\beta_n|^r<\oo$.  If\/ $\sum_{n=1}^\oo \tau_n \beta_n<\oo$, then
$\sum_{n=1}^\oo \tau_n \alpha_n^r < \infty$.
\end{lem}

\begin{proof}
There is a polynomial
$p_r(x,y)$ of degree $r-1$ with coefficients depending only on $r$ such that
$(x-y)^r=x^r-yp_r(x,y)$.  Let $c_r$ be the maximum of $p_r$ over
$[0,1]^2$.  Then:
$$
   \oo > \sum_{n=1}^\oo \tau_n |\alpha_n-\beta_n|^r
    =\sum_{n=1}^\oo \tau_n |\alpha_n^r-\beta_n p_r(\alpha_n,\beta_n)|
    \ge \sum_{n=1}^\oo \tau_n (\alpha_n^r - c_r \beta_n).
$$
If $\sum_{n=1}^\oo \tau_n \beta_n$ converges, then it follows that
$\sum_{n=1}^\oo \tau_n \alpha_n^r$ also converges.
\end{proof}

\begin{lem}\label{lem:the-lemma}
  Let $\{ \tau_n \}$ be a positive sequence and let
  $\{ b_n \}$ be a strictly positive increasing sequence for $n\ge 1$.
  Let $X$ be a random variable
  such that
\begin{equation}\label{eq:crit}
  \sum_{n=1}^\oo \tau_n n P(|X|\ge b_n) < \oo.
\end{equation}
  Fix $\nu\in[0,\oo)$.	Put $T_k=\sum_{n=1}^k n\tau_n$.
  Suppose that there is a constant $C\in (0,\oo)$ and a $\theta \in
  [1,\oo)$ such that:
\begin{equation}\label{eq:aux-gen-1}
  \frac{b_n^{\nu\theta}}{n^{\theta-1}}\sum_{k=n}^\oo \frac{k^\theta \tau_k}{b_k^{\nu\theta}}
    \le C T_{n-1},
\end{equation}
for all $n\ge 2$, and
\begin{equation}\label{eq:aux-gen-2}
   \frac{k b_n^\nu}{b_k^\nu} \le C T_{n-1},
\end{equation}
whenever $k\ge n\ge 2$.  Then:
$$
\sum_{n=1}^\oo \tau_n \left(\frac{nE[|X|^\nu\cdot 1_{\{ |X_1|<b_n \}}]}
	       {b_n^\nu}\right)^{\!\!\theta}<\oo
$$
\end{lem}

The proof is based on methods of Klesov~\cite[Proof of Theorem~4]{Klesov}.

\begin{proof}
Set $X^{(n)}=X\cdot 1_{\{ |X|<b_n \}}$.
Let $b_0=0$ and put $t_n=E[|X^{(n)}|^\nu]$ for $n\ge 0$.  Note that $t_0=0$.
Put $\delta_n=t_n^\theta-t_{n-1}^\theta$ for $n\ge 1$.
For convenience, let $T_0=1$;  redefining $C$ if necessary, we may assume
that \eqref{eq:aux-gen-1} and \eqref{eq:aux-gen-2} also hold for $n=1$.
Then, since $t_n^\theta=\sum_{k=1}^n
\delta_k$, and using Fubini's theorem and \eqref{eq:aux-gen-1}:
\begin{equation}\label{eq:Klesov-1}
\begin{split}
\sum_{n=1}^\oo \tau_n
	\left(\frac{nE[|X^{(n)}|^\nu]}{b_n^\nu}\right)^{\!\!\theta}
      &= \sum_{n=1}^\oo \tau_n (nb_n^{-\nu})^\theta \sum_{k=1}^n \delta_k \\
      &= \sum_{k=1}^\oo \delta_k \sum_{n=k}^\oo \tau_n (nb_n^{-\nu})^\theta \\
      &\le C\sum_{k=1}^\oo \delta_k k^{\theta-1} b_k^{-\nu\theta} T_{k-1}.
\end{split}
\end{equation}
Let $P_k=P(b_{k-1}\le |X| < b_k)$.  Then, note that
\begin{equation}\label{eq:delta-k}
	\delta_k \le c (t_k-t_{k-1})t_k^{\theta-1} \le cP_k b_k^\nu t_k^{\theta-1},
\end{equation}
where $c$ is a finite constant depending only on $\theta$ and such
that $x^\theta-y^\theta\le c (x-y)x^{\theta-1}$ whenever $0\le y\le x$.
Now, fix $\ell\in\Z^+$.  Let $\rho_n=\ell/b_\ell^\nu$ for $n=\ell$ and put $\rho_n=0$ for all other
$n$.  Observe that by \eqref{eq:aux-gen-2} we have
\begin{equation}\label{eq:check-lelem}
	b_n^\nu \sum_{j=n}^\oo \rho_j  \le C T_{n-1},
\end{equation}
for $n\le \ell$ and the same inequality trivially holds for $n>\ell$.

By Lemma~\ref{lem:elementary} (with $t=\nu$) and \eqref{eq:check-lelem},
it follows that if $\ell \ge 2$ so that $\rho_1=0$, then
then
\[
\begin{split}
	\frac{lt_l}{b_l^\nu}
	&=\sum_{n=1}^\oo \rho_n t_n
	=\sum_{n=2}^\oo \rho_n E[|X^{(n)}|^\nu] \\
	&\le C \sum_{n=1}^\oo n\tau_n
	P(|X|\ge b_n)
	<\oo,
\end{split}
\]
where the finiteness of the right hand side followed from \eqref{eq:crit}.
Thus, $K\eqdef\sup_{\ell\ge 1} \ell t_\ell/b_\ell^\nu<\oo$.
Then, $t_k \le K
b_k^\nu/k$, so that \eqref{eq:delta-k} yields:
$$
	\delta_k \le cK^{\theta-1} P_k b_k^\nu (b_k^\nu/k)^{\theta-1}=cK^{\theta-1} \frac{P_k b_k^{\nu\theta}}{k^{\theta-1}}.
$$
Putting this into \eqref{eq:Klesov-1}, and recalling that $T_1=1$, we see that:
\[
\begin{split}
	\sum_{n=1}^\oo \tau_n
	\left(\frac{nE[|X^{(n)}|^\nu]}{b_n^\nu}\right)^{\!\!\theta}
      &\le CcK^{\theta-1}\sum_{k=1}^\oo P_k T_{k-1} \\
      &=CcK^{\theta-1}\left(P_1+\sum_{k=2}^\oo P_k \sum_{n=1}^{k-1} n\tau_n)\right) \\
      &=CcK^{\theta-1}\left(P_1+\sum_{n=1}^\oo n\tau_n \sum_{k=n+1}^\oo P_k\right) \\
      &=CcK^{\theta-1}\left(P_1+\sum_{n=1}^\oo n\tau_n P(|X|\ge b_n)\right) < \oo,
\end{split}
\]
by Fubini's theorem and \eqref{eq:crit}.
\end{proof}

The following version of the Hoffman-J\o rgensen inequality~\cite{HJ} will be
needed for the proof of Theorem~\ref{th:main} in the case $\theta>1$ and
follows immediately from \cite[Lemma~2.2]{LRJW}.

\begin{lem}\label{lem:HJ} Let $X_1,\dots,X_n$ be independent symmetric random
variables, and let $S_n=X_1+\dots+X_n$.  Then for each $r \in\Z^+$ there
exist finite constants $C_r$ and $D_r$ such
that for all $\lambda\ge 0$ we have:
$$
	P(|S_n|\ge \lambda)
	 \le C_r \sum_{k=1}^n P(|X_k|\ge \lambda/(2r))
	 + D_r [P(|S_n|\ge \lambda/(2r))]^r.
$$
\end{lem}

Let $\Phi$ be the distribution function of a (0,1) normal random variable.

\begin{lem}\label{lem:2}
Under the global conditions of Theorem~\upref{th:main},
suppose that condition \upn{(ii)} is satisfied and that
$X_1$ is symmetric.
Then, the following four conditions are equivalent:
\begin{enumerate}
\item[\upn{(a)}] $\sum_{n=1}^\oo \tau_n P(|S_n|\ge \e a_n)<\oo$ for
all $\e>0$;
\item[\upn{(b)}] for all $\gamma \in (0,1]$ and $\e>0$ we have
$\sum_{n=1}^\oo \tau_n [1-\Phi(\gamma\e a_n/(nT_{\e,n})^{1/2})]<\oo$,
where $T_{\e,n}$ is as in Theorem~\upref{th:main},
and where $\Phi(t/0)=0$ for all $t>0$;
\item[\upn{(c)}] there is an $s\ge 1$ such that for all $\e>0$ we have
$\sum_{n=1}^\oo \tau_n [1-\Phi(\e a_n/(nT_{\e,n})^{1/2})]^s<\oo$;
\item[\upn{(d)}] condition~\upn{(iii)} of Theorem~\upref{th:main} holds.
\end{enumerate}
\end{lem}

The equivalence of (a) and (b) will be the most difficult part to prove,
and will involve a similar method of proof to that in~\cite{Pruss:Riemann},
using a central limit theorem estimate in the present case due to
Nagaev~\cite{Nagaev}.  In the case where $\theta>1$, we will also need
the Hoffman-J\o rgensen inequality based methods of Klesov~\cite[Proof
of Theorem~4]{Klesov}.

\begin{proof}[Proof of Lemma~\upref{lem:2}]
Note that:
\begin{equation}\label{eq:Phi-ineq}
	1-\Phi(x)=\pi^{-1/2} x^{-1} e^{-x^2/2} (1-O(x^{-2}))
\end{equation}
as $x\to\oo$.  Thus:
\begin{equation}\label{eq:Phi-ineq2}
	1-\Phi(x)\le C e^{-x^2/2},
\end{equation}
for all $x\ge 0$ and an absolute constant $C$.

First suppose (d) holds.
Put $\gamma'=\gamma/2$.  Then, (d) implies that:
\begin{equation}\label{eq:ge2}
	\sum_{n=1}^\oo \tau_n e^{-(\gamma'\e)^2 a_n^2
	/nT_{\gamma'\e,n}}<\oo,
\end{equation}
for all $\e>0$.  Observe now that as $\gamma'<\gamma\le 1$ and by
definition of $T_{\e,n}$:
\begin{equation}\label{eq:Tsplit}
\begin{split}
	T_{\e,n} &\le T_{\e\gamma',n} + (\e\gamma a_n)^2
	P(|X_1|\ge \e\gamma' a_n) \\
	&\le \max(2T_{\e\gamma',n},2(\e \gamma a_n)^2 P(|X_1|\ge
	\e\gamma' a_n)).
\end{split}
\end{equation}
Then, by \eqref{eq:Phi-ineq2} and \eqref{eq:Tsplit}:
\[
\begin{split}
	\sum_{n=1}^\oo \tau_n &[1-\Phi(\gamma\e a_n / (nT_{\e,n})^{1/2})]
	  \le C \sum_{n=1}^\oo \tau_n e^{-\gamma^2\e^2a_n^2/2nT_{\e,n}} \\
	  &\le C\max\left( \sum_{n=1}^\oo \tau_n
	  e^{-\gamma^2\e^2a_n^2/4nT_{\e\gamma',n}}
	  , \sum_{n=1}^\oo \tau_n e^{-1 /
	  4 n P(|X_1|\ge\e\gamma' a_n)} \right),
\end{split}
\]
where we use the convention that $e^{-1/0} = 0$.
Now, the first sum on the right hand side here converges by
\eqref{eq:ge2} since $\gamma^2/4=(\gamma')^2$, while
the second converges by
the elementary inequality $e^{-1/x} \le x$ valid for all
$x\ge 0$, together with condition (ii) of Theorem~\upref{th:main}.
Hence, (d) implies (b).

The implication from (b) to (c) is trivial.  Suppose now that (c) holds.
Let $u_{\delta,n}=a_n/(nT_{\delta,n})^{1/2}$.
Then,
\begin{equation}\label{eq:from-b}
	\sum_{n=1}^\oo \tau_n
	 e^{-s \delta^2 u_{\delta,n}^2} < \oo,
\end{equation}
for all $\delta>0$ by (c) and \eqref{eq:Phi-ineq}.  Fix $\e>0$ and let
$\delta=\e/(2s)^{1/2}$.  Then, $u_{\delta,n}\le u_{\e,n}$ and so
$$
e^{-s \delta^2 u_{\delta,n}^2} \ge e^{-s \delta^2 u_{\e,n}^2 }
 = e^{-\e^2 u_{\e,n}^2 / 2 }.
$$
Condition (d) then follows from this together with
\eqref{eq:Phi-ineq} and \eqref{eq:from-b}.
%
%
%
Hence (c) implies (d), so that we have shown that
$\text{(d)}\Rightarrow\text{(b)}\Rightarrow\text{(c)}\Rightarrow\text{(d)}$.

All we now need to prove is the equivalence of (a) and (b).
To do this, assume we are in the setting
of Theorem~\ref{th:main} and that (ii) holds.
Fix $\e>0$.  Let
$X_{k,n}(\e)=X_{k}\cdot 1_{\{ |X_k|<\e a_n \}}$.  Let
$S_n(\e)=X_{1,n}(\e)+\dots+X_{n,n}(\e)$.  Set $A_n(\e)=\bigcup_{k=1}^n \{
X_k\ne X_{k,n}(\e) \}$.
Observe that $S_n=S_n(\e)$ except possibly on $A_n(\e)$, and that
\begin{equation}\label{eq:A-sum}
	\sum_{n=1}^\oo \tau_n P(A_n(\e)) \le \sum_{n=1}^\oo n\tau_n
	P(|X_1|\ge \e a_n) < \oo,
\end{equation}
by (ii).  It follows from \eqref{eq:A-sum} and from the equality of
$S_n$ and $S_n(\e)$ outside $A_n(\e)$ that (a) holds if and only if
\begin{equation}\label{eq:convp}
	\sum_{n=1}^\oo \tau_n P(|S_n(\e)|\ge \e a_n)<\oo,\qquad\forall
	\e>0.
\end{equation}

We now need only show that \eqref{eq:convp} holds if and only if (b)
holds, and we will be done.  Note that $T_{\e,n}=E[(X_{1,n}(\e))^2]$.
Fix $\gamma>0$ to be chosen later as needed.  Then,
since all the $X_{k,n}(\e)$ and $S_n(\e)$ have mean zero by symmetry, and since
$X_{1,n}(\e),\dots,X_{n,n}(\e)$ are identically distributed for a fixed $n$, by
Nagaev's central limit theorem estimate~\cite[Theorem~3]{Nagaev} we have:
$$
	\left|P(|S_n(\e)|\ge \gamma\e a_n)-2\left[1-\Phi\left(
\frac{\gamma\e a_n}{nT_{\e,n}^{1/2}}\right)\right]\right|
	 \le c n \frac{E[|X_{1,n}(\e)|^3]}{(\gamma\e a_n)^3},
$$
for an absolute constant $c<\oo$.  Since the left hand side never exceeds $1$,
it follows further that for a (possibly different) absolute constant $c$, we have:
\begin{equation}\label{eq:Nagaev}
	\left|P(|S_n(\e)|\ge \gamma\e a_n)-2\left[1-\Phi\left(
\frac{\gamma\e a_n}{nT_{\e,n}^{1/2}}\right)\right]\right|
         \le c \min\left(1, n \frac{E[|X_{1,n}(\e)|^3]}{(\gamma\e a_n)^3}
	 \right),
\end{equation}

Observe that if \eqref{eq:aux-cond} holds for some $\theta$, then
it also holds for all greater values.
We now have a quick proof if $\theta=1$.  For then, by \eqref{eq:Nagaev}
(with $\gamma=1$)
we
can see that the equivalence of (b) and \eqref{eq:convp} would follow
as soon as we could show that we have
\begin{equation}\label{eq:final-conv}
	\sum_{n=1}^\oo n \tau_n \frac{E[|X_{1,n}(\e)|^3]}{a_n^3}<\oo.
\end{equation}
But \eqref{eq:final-conv} follows from the validity of
\eqref{eq:aux-cond} for $\theta=1$ and from condition~(ii) of
Theorem~\ref{th:main}, by an application of Lemma~\ref{lem:elementary}
with $\rho_n=n\tau_n/a_n^3$, $X\equiv X_1$, $t=3$ and  $b_n=\e a_n$.

Suppose now we are working with $\theta>1$, so that
\eqref{eq:aux-aux-cond} also holds.  Let $r$ be an integer greater than
or equal to $\theta$.  I now claim that:
\begin{equation}\label{eq:the-claim}
	\sum_{n=1}^\oo \tau_n
        \left(\min\left(1,\frac{nE[|X_{1,n}(\e)|^3]}{a_n^3}\right)\right)^{\!\!r} < \oo, \qquad
	\forall \e>0.
\end{equation}

Suppose for now that this has been shown.  If (a) holds, then as noted
before, \eqref{eq:convp} does likewise.  Letting
$\alpha_n=2[1-\Phi(\e a_n/(nT_{\e,n})^{1/2}))]$ and
$\beta_n=P(|S_n(\e)|\ge \e a_n)$, we see that by \eqref{eq:convp}
together with \eqref{eq:the-claim} and \eqref{eq:Nagaev} (with
$\gamma=1$), we do have the conditions of Lemma~\ref{lem:comp}
satisfied, so that $\sum_{n=1}^\oo \tau_n \alpha_n^r < \oo$, and hence
(c) follows, whence (b) follows by the already proved equivalence of
(b), (c) and (d).

Conversely, suppose (b) holds.	This time letting
$\beta_n=2[1-\Phi(\gamma\e a_n/(nT_{\e,n})^{1/2}))]$ and
$\alpha_n=P(|S_n(\e)|\ge \gamma\e a_n)$, using (b), together with
\eqref{eq:Nagaev}, \eqref{eq:the-claim} and Lemma~\ref{lem:comp}, we see
that $\sum_{n=1}^\oo \tau_n \alpha_n^r<\oo$, i.e.,
\begin{equation}\label{eq:res1}
   \sum_{n=1}^\oo \tau_n [P(|S_n(\e)|\ge \gamma\e a_n)]^r<\oo.
\end{equation}
Let $\gamma=(2r)^{-1}$.  By the Hoffman-J\o rgensen inequality
(Lemma~\ref{lem:HJ}), we have:
\[
\begin{split}
  \sum_{n=1}^\oo &\tau_n P(|S_n(\e)|\ge \e a_n) \\
    &\le C_r \sum_{n=1}^\oo \tau_n P(|X_n(\e)|\ge \gamma\e a_n)
     + D_r \sum_{n=1}^\oo \tau_n [P(|S_n(\e)|\ge \gamma\e a_n)]^r.
\end{split}
\]
It is easy to see that the first sum on the right hand side is no greater than $\sum_{n=1}^\oo
\tau_n P(A_n(\gamma\e))$, which converges by \eqref{eq:A-sum}, and the
second converges by \eqref{eq:res1}.  Hence, \eqref{eq:convp} follows,
and as already shown this implies (a).

Hence, all we need to show is that \eqref{eq:the-claim} holds.
Changing finitely many values of $a_n$ and $\tau_n$, we may assume that
\eqref{eq:aux-cond} holds for all $n\ge 2$, and that \eqref{eq:aux-aux-cond}
gives:
\begin{equation}\label{eq:aux-cond-mod}
	\sup_{k\ge n} \frac{k a_n^3}{a_k^3}
	\le C'\sum_{k=1}^{n-1} n\tau_n,
\end{equation}
for all $n\ge 2$ and a finite $C'$.  We can now apply Lemma~\ref{lem:the-lemma} with $\nu=3$,
$b_n=\e a_n$, and $X\equiv X_1$, using the assumed condition~(ii) of
Theorem~\ref{th:main} to guarantee \eqref{eq:crit}, and getting
\eqref{eq:aux-gen-1} and \eqref{eq:aux-gen-2} (with an appropriately
chosen constant) from \eqref{eq:aux-cond} and
\eqref{eq:aux-cond-mod} (which hold for all $n\ge 2$ by
assumption),
respectively.  The Lemma then yields \eqref{eq:the-claim} since $r\ge\theta>1$,
and
the proof is complete.
\end{proof}

\begin{proof}[Proof of Theorem~\upref{th:main}]
Let $\alpha_1,\alpha_2,\dotsc$ be a sequence of independent Bernoulli
random variables with $P(\alpha_n=1)=P(\alpha_n=-1)=\half$, and with the
sequence independent of $\{ X_n \}_{n=1}^\oo$. Let $X_n'=\alpha_n X_n$,
and put $S_n'=X_1'+\dots+X_n'$.
Note that the primed variables are symmetric.

Suppose first that \eqref{eq:conv} holds and $\{ \tau_n \}$ satisfies
Condition~A.  By Lemma~\ref{lem:1} we have (i) and (ii) holding.  It
remains to show that (iii) holds.    But, if \eqref{eq:conv} holds, it
likewise holds with $S_n'$ in place of $S_n$, as can be seen by
conditioning on the $\{ \alpha_k \}_{k=1}^\oo$ and
using~\cite[Corollary~5]{MS}. But then by the implication
$\text{(a)}\Rightarrow\text{(d)}$ of Lemma~\ref{lem:2} applied to the
primed variables (which are symmetric) it follows that (iii) holds with
$X_1'$ in place of $X_1$.  But (iii) holding for $X_1'$ is equivalent to
it holding for $X_1$, and so (iii) follows.

Conversely, suppose (i), (ii) and (iii) hold.  Evidently (ii) and (iii)
will also
hold with $X_1'$ in place of $X_1$.  Hence by
the implication $\text{(d)}\Rightarrow\text{(a)}$ of Lemma~\ref{lem:2}
as applied to the primed variables we have
\begin{equation}\label{eq:symm-conv}
	\sum_{n=1}^\oo \tau_n P(|S_n'|\ge \e a_n)<\oo, \qquad
	\forall\e>0.
\end{equation}
Now, for any random variable $Y$, let $Y^s=Y-\tilde Y$ be the
symmetrization of $Y$, where $\tilde Y$ is an independent copy of $Y$.
We shall choose symmetrizations in such a way that $X_1^s,X_2^s,\dotsc$
are independent and $S_n^s=X_1^s+\dots+X_n^s$ for all $n$.  Now,
since
\[
\begin{split}
	P(|X_k^s|\ge \lambda)&\le P(|X_k|\ge \lambda/2)+P(|\tilde
	X_k|\ge \lambda/2)\\
	&=2P(|2X_k|\ge \lambda)
	=2P(|2X_k'|\ge \lambda),
\end{split}
\]
for all $\lambda\ge 0$, it follows from \cite[Theorem~1]{Pruss:noniid}
that there is an absolute constant $c>0$ such that
$$
	P(|S_n^s|\ge \lambda)\le c P(|2S_n'|\ge \lambda/c)
$$
for all $\lambda\ge 0$.  By \eqref{eq:symm-conv} it then follows that
$$
	\sum_{n=1}^\oo \tau_n P(|S_n^s|\ge \e a_n)<\oo,
	\qquad\forall\e>0.
$$
By standard symmetrization inequalities~\cite[\S17.1.A]{Loeve} it follows
that
\begin{equation}\label{eq:symm-conv2}
	\sum_{n=1}^\oo \tau_n P(|S_n-\mu_n|\ge \e a_n)<\oo,
	\qquad\forall\e>0,
\end{equation}
where $\mu_n$ is the median of $S_n$ occurring in (i).
Now, for any $\e>0$ we have
\begin{multline}\label{eq:for-ref}
	\sum_{n=1}^\oo \tau_n P(|S_n|\ge \e a_n) \\
	 \le \sum_{n=1}^\oo \tau_n P(|S_n-\mu_n|\ge \e a_n/2)
	   + \sum_{n=1}^\oo \tau_n \cdot 1_{\{ |\mu_n|>\e a_n/2 \}}.
\end{multline}
By condition~(i), the second sum converges, and by
\eqref{eq:symm-conv2}, so does the first, and hence \eqref{eq:conv}
follows.
\end{proof}

\begin{proof}[Proof of Theorem~\upref{th:aux}]
Assume condition~(ii) of Theorem~\ref{th:main} holds.  Fix $\e>0$.
Changing a finite number of values of $\tau_n$ and $a_n$
and using \eqref{eq:aux-cond2} and \eqref{eq:aux-aux-cond2}
will let us assume that if we let $b_n=\e a_n$ and $\nu=2$, then
conditions \eqref{eq:aux-gen-1} and
\eqref{eq:aux-gen-2} of Lemma~\ref{lem:the-lemma} will be verified
for an appropriate constant.
Applying that Lemma and using condition~(ii) of
Theorem~\ref{th:main} shows that:
\begin{equation}\label{eq:t-aux-1}
	\sum_{n=1}^\oo \tau_n n a_n^{-2} T_{\e,n} < \oo,
\end{equation}
where $T_{\e,n}$ is as in Theorem~\ref{th:main}(iii).
Now, using the elementary inequality $e^{-1/x} \le c_\theta x^\theta$ which
is valid for all $x\ge 0$ where $c_\theta$ is a constant depending only on
$\theta>0$, we see that
\eqref{eq:t-aux-1} entails that
$$
\sum_{n=1}^\oo \tau_n e^{-\e^2 a_n^2 / (nT_{\e,n})} < \oo,
$$
and so (iii) is true.
\end{proof}

\begin{proof}[Proof of Lemma~\upref{lem:Sp}]
Let $a_n=(n\log n)^{1/2}$ for $n\ge 2$.  Assume the conditions of the
Lemma hold.
Without loss of generality $E[X_1]=0$.
Let $Y_n=X_1\cdot 1_{\{|X_1|\le a_n\}}$.
Observe that Condition~(b) of Corollary~\ref{cor:Sp} implies that
\begin{equation}\label{eq:show-1}
\sum_{n=2}^\oo P(|X_1|> a_n)<\oo.
\end{equation}

We now need the fact that
in order to prove that $S_n/a_n\to 0$ in probability, given the fact
that the $X_k$ are identically distributed and with mean zero, all we
need to show is that \eqref{eq:show-1} holds as well as that
\begin{equation}\label{eq:show-2}
	\lim_{n\to\oo} \frac{n E[Y_n^2]}{a_n^2} \to 0
\end{equation}
and
\begin{equation}\label{eq:show-3}
	\lim_{n\to\oo} \frac{n E[Y_n]}{a_n} \to 0.
\end{equation}
To see this fact, let $\xi_{n,k}=X_k/a_n$ and $m_{n,k}=0$ in the setting of
\cite[p.~105]{GnedenkoKolmogorov} (there $m_{n,k}$ is supposed to be a median
of $\xi_{n,k}$, but the proof in the direction we need does not use this
fact), and note that \eqref{eq:show-1}--\eqref{eq:show-3} together with
the assumption that $E[X_1]=0$ (which assumption implies that $E[X_1\cdot
1_{\{|X_1| > 1}]=-E[Y_1]$) and the proof
in~\cite[p.~106]{GnedenkoKolmogorov} yield the convergence in
probability of $(S_n-nE[Y_n])/a_n$ to zero, which in turn implies the
convergence in probability of $S_n/a_n$ to zero by \eqref{eq:show-3}.

Since as already noted, Condition~(b) of Corollary~\ref{cor:Sp} implies
\eqref{eq:show-1}, we need only
prove \eqref{eq:show-2} and \eqref{eq:show-3}.
First, observe that for $n\ge 2$ we have:
\begin{equation}\label{eq:showing-2}
\begin{split}
	\frac{n E[Y_n^2]}{a_n^2}
	 &= (\log n)^{-1} E[X_1^2 \cdot 1_{\{|X_1|\le a_n\}}] \\
	 &=E\left[\frac{X_1^2}{\log^+ X_1}\cdot
	 \frac{1_{\{ |X_1|\le a_n\}} \cdot \log^+X_1}{\log n} \right] .
\end{split}
\end{equation}
Let
$$
	W_n=\frac{X_1^2}{\log^+ X_1}\cdot
         \frac{ 1_{\{ |X_1|\le a_n\}} \cdot \log^+X_1}{\log n}.
$$
Observe that
$$
	0\le W_n \le c \frac{X_1^2}{\log^+ X_1},
$$
where $c=\sup_{n\ge 2} (\log^+ a_n)/(\log n)<\oo$.  But $X_1^2/\log^+
X_1$ is in $L^1$, while $W_n\to 0$ pointwise as $n\to\oo$, so that by
dominated convergence it follows that $E[W_n]\to 0$ as $n\to\oo$ and
hence \eqref{eq:show-2} follows from \eqref{eq:showing-2}.

Now, to prove \eqref{eq:show-3}, let $U_n=X_1-Y_n$.
Observe
that $E[U_n]=-E[Y_n]$ since $E[X_1]=0$.  Hence,
we only need to show that $a_n^{-1} n
E[U_n]\to 0$.  Now, since $|U_n|\ge a_n$ whenever $U_n\ne 0$
and as $x^{-1}\log^+ x$ is decreasing in $x>0$, we have
(using the convention that $(0^2/\log^+ |0|)\cdot (\log^+ |0|)/|0|=0$):
\[
\begin{split}
	a_n^{-1} n E[|U_n|]
	  &= a_n^{-1} n E[(U_n^2 / \log^+ |U_n|) \cdot (\log^+ |U_n|) /
	  |U_n|]
	  \\
	  &\le a_n^{-1} n E[(U_n^2/\log^+ |U_n|)\cdot (\log^+ a_n)/a_n] \\
	  &= (\log^+ a_n) n a_n^{-2} E[(U_n^2/\log^+ U_n)] \to 0,
\end{split}
\]
since $(\log^+ a_n)na_n^{-2}$ is bounded while $E[(U_n^2/\log^+ U_n)]\to
0$ by dominated convergence, as $U_n\to 0$ almost surely and
$U_n^2/\log^+ U_n \le X_1^2/\log^+ X_1$, while $X_1^2/\log^+ X_1 \in L^1$ by
Condition~(b) of Corollary~\ref{cor:Sp}.
\end{proof}

\section*{Acknowledgments}
The author is most grateful for a number of interesting e-mail
conversations on these topics with Professors Aurel Sp\u ataru and
Stephen Montgomery-Smith.  In particular, the author is very grateful to
Professor Montgomery-Smith for permission to use his
Example~\ref{ex:MS}, and to Professor Sp\u ataru for a number of useful
comments on earlier versions of this paper, including pointing out
several errors.  The author also expresses his gratitude for a very careful
reading to two referees, and thanks for pointing out a number of errors in
earlier drafts.

\providecommand{\bysame}{\leavevmode\hbox to3em{\hrulefill}\thinspace}

\end{document}